\documentclass[12pt]{article}

\usepackage{amsmath}
\usepackage{amssymb}
\usepackage{amsthm}
\usepackage{bbm}
\usepackage{graphicx,psfrag}
\usepackage{enumerate}
\usepackage{url} 
\usepackage{bm}
\usepackage{dsfont}
\usepackage{booktabs}
\usepackage{caption}
\usepackage{hyperref}
\usepackage{amssymb}
\usepackage{wasysym}
\usepackage{color}
\usepackage{mathtools}
\newcommand{\bracks}[1]{\left[#1\right]}
 
\newcommand{\parens}[1]{\left(#1\right)}

\newcommand{\R}{{\mathbb R}}

\newcommand{\bbR}{\mathbb{R}}
\newcommand{\bbM}{\mathbb{M}}
\newcommand{\bbS}{\mathbb{S}}
\newcommand{\sdp}{ \bbS^d_+}

\newcommand{\trace}{\mathrm{tr\,}}

\newcommand{\csch}{\mathrm{csch}\,}
\newcommand{\define}{\coloneqq}

\newcommand{\BM}{\mathrm{Bm}}
\newcommand{\CBM}{\mathrm{cBm}}
\newcommand{\OU}{\mathrm{OU}}

\newcommand{\BB}{\mathrm{Bb}}

\newcommand{\half}{ {\scriptstyle{\frac{1}{2}} } }
\newcommand{\Var}{\mathrm{Var}}
\newcommand{\beqn}{\vspace{-0.25cm}\begin{eqnarray*}}
\newcommand{\eeqn}{\end{eqnarray*}}
\newcommand{\bneqn}{\vspace{-0.25cm}\begin{eqnarray}}
\newcommand{\eneqn}{\end{eqnarray}}
\newcommand{\convd}{~{\buildrel \mathcal{D} \over \rightarrow}~}
\newcommand{\convas}{~{\buildrel a.s. \over \longrightarrow}~}

\newtheorem{proposition}{Proposition}
\newtheorem{theorem}{Theorem}

\theoremstyle{definition}
\newtheorem{remark}{Remark}

\def\tcp{\textcolor{black}}

\title{Yule's ``nonsense correlation'' solved: Part II}
\author{Philip A. Ernst\footnote{Department of Mathematics, Imperial College London. Email:\,p.ernst@imperial.ac.uk}\,\,, L.C.G. Rogers\footnote{Statistical Laboratory, University of Cambridge. Email: chris@statslab.cam.ac.uk }\,\,, and Quan Zhou\footnote{Department of Statistics, Texas A\&M University. Email: quan@stat.tamu.edu}}
\date{\today}
\begin{document}
\maketitle

\begin{abstract}
\tcp{In 1926, G. Udny Yule (\cite{yule1926}) considered the following problem: given a sequence of pairs of random variables $\{X_k,Y_k \}$ ($k=1,2, \ldots, n$), and letting $X_i = S_i$ and $Y_ i= S'_i$ where $S_i$ and $S'_i$ are the partial sums of two independent random walks, what is the distribution of the empirical correlation coefficient
\begin{equation*}
\rho_n = \frac{\sum_{i=1}^n S_i S^\prime_i - \frac{1}{n}(\sum_{i=1}^n S_i)(\sum_{i=1}^n S^\prime_i)}{\sqrt{\sum_{i=1}^n S^2_i - \frac{1}{n}(\sum_{i=1}^n S_i)^2}\sqrt{\sum_{i=1}^n (S^\prime_i)^2 - \frac{1}{n}(\sum_{i=1}^n S^\prime_i)^2}}?
\end{equation*}
Yule empirically observed the distribution of this statistic to be heavily dispersed and frequently large in absolute value, leading him to call it ``nonsense correlation.'' This unexpected finding led to his formulation of two concrete questions, each of which would remain open for more than ninety years: (i) Find (analytically) the variance of $\rho_n$ as $n \rightarrow \infty$ and (ii): Find (analytically) the higher order moments and the density of $\rho_n$ as $n \rightarrow \infty$. In 2017, the authors of \cite{ernst2017yule} considered the empirical correlation coefficient
\begin{equation*}
\rho:= \frac{\int_0^1W_1(t)W_2(t) dt - \int_0^1W_1(t) dt \int_0^1 W_2(t) dt}{\sqrt{\int_0^1 W^2_1(t) dt - \parens{\int_0^1W_1(t) dt}^2} \sqrt{\int_0^1 W^2_2(t) dt - \parens{\int_0^1W_2(t) dt}^2}}
\end{equation*}
of two \textit{independent} Wiener processes $W_1,W_2$, the limit to which $\rho_n$ converges weakly, as was first shown by the author of \cite{phil}. Using tools from integral equation theory, the authors of \cite{ernst2017yule} closed question (i) by explicitly calculating the second moment of $\rho$ to be .240522. This paper begins where \cite{ernst2017yule} leaves off. We succeed in closing question (ii) by explicitly calculating all moments  of $\rho$ (up to order 16).}
This leads, for the first time, to an approximation to the density of Yule's nonsense correlation. We then proceed  explicitly to compute higher moments of  $\rho$ when the two independent Wiener processes are replaced by two correlated Wiener processes, two independent Ornstein-Uhlenbeck processes, and two independent Brownian bridges. 
We conclude by extending the definition of $\rho$ to the time interval $[0, T]$ for any $T > 0$ and prove a Central Limit Theorem for the case of two independent Ornstein-Uhlenbeck processes.
\end{abstract}

\section{Introduction}
\tcp{Given a sequence of pairs of random variables $\{X_k,Y_k \}$ ($k=1,2, \ldots, n$), how can we measure the strength of the dependence of $X$ and $Y$? } The classical Pearson correlation coefficient addresses this question in
the case when the sequence $(X_k, Y_k)_{k = 1}^n$ is an IID sequence, but
when blindly calculated for sequences which are not IID, as arise (for
example) when the sequence $(X_k, Y_k)_{k = 1}^n$ is a random walk or other
time series, the resulting statistic is meaningless. \tcp{Letting $X_i = S_i$ and $Y_ i= S'_i$, where $S_i$ and $S'_i$ are the partial sums of two independent random walks, Yule \cite{yule1926} considered the distribution of the empirical correlation coefficient
\begin{equation*}
\rho_n = \frac{\sum_{i=1}^n S_i S^\prime_i - \frac{1}{n}(\sum_{i=1}^n S_i)(\sum_{i=1}^n S^\prime_i)}{\sqrt{\sum_{i=1}^n S^2_i - \frac{1}{n}(\sum_{i=1}^n S_i)^2}\sqrt{\sum_{i=1}^n (S^\prime_i)^2 - \frac{1}{n}(\sum_{i=1}^n S^\prime_i)^2}},
\end{equation*}
which he found to be both heavily dispersed and frequently large in absolute value. This led Yule to call this
distribution ``nonsense correlation''\footnote{In lieu of calling this distribution ``nonsense correlation,'' the authors of \cite{ernst2017yule} refer to it as ``volatile'' correlation since its distribution is both heavily dispersed and is frequently large in absolute value.} and to formulate two concrete questions about its distribution, each of which would remain open for more than ninety years: 
\begin{enumerate}[(i)]
\item Find (analytically) the variance of $\rho_n$ as $n \rightarrow \infty$.
\item Find (analytically) the higher order moments and the density of $\rho_n$ as $n \rightarrow \infty$.
\end{enumerate} 
}

\indent Despite the prominence of Yule's 1926 paper, his findings would remain ``isolated'' from the literature until 1986 (see \cite{Aldrich}), when the authors of \cite{Hendry} and \cite{phil} confirmed many of the empirical claims of ``spurious regression'' made by the authors of \cite{gran}. In particular, \cite{phil} provided a mathematical solution to the problem of spurious regression among integrated time series by demonstrating that statistical {\em t}-ratio and {\em F}-ratio tests diverge with the sample size, thereby explaining the observed `statistical significance' in such regressions. In later work \cite{phil2}, the same author provided an explanation of such spurious regressions in terms of orthonormal representations of the Karhunen Lo\`{e}ve type.\\ 
\tcp{\indent In 2017, the authors of \cite{ernst2017yule} considered the empirical correlation coefficient
\begin{equation*}
\rho:= \frac{\int_0^1W_1(t)W_2(t) dt - \int_0^1W_1(t) dt \int_0^1 W_2(t) dt}{\sqrt{\int_0^1 W^2_1(t) dt - \parens{\int_0^1W_1(t) dt}^2} \sqrt{\int_0^1 W^2_2(t) dt - \parens{\int_0^1W_2(t) dt}^2}}
\end{equation*}
of two \textit{independent} Wiener processes $W_1,W_2$, the limit to which $\rho_n$ converges weakly, as originally shown in \cite{phil}. Using tools from integral equation theory, \cite{ernst2017yule} closed question (i) by explicitly calculating the second moment of $\rho$ to be .240522.\\}
\indent \tcp{The present paper begins where \cite{ernst2017yule} leaves off by closing question (ii).} We remarkably succeed in closing this longstanding open question by explicitly calculating all moments  of $\rho$ (up to order 16). These lead, for the first time, to an approximation to the density of Yule's nonsense correlation. We then proceed to go further beyond Yule's original message by explicitly computing higher moments of the empirical correlation coefficient when the two independent Wiener processes are replaced by two correlated Wiener processes, two independent Ornstein-Uhlenbeck processes, and two independent Brownian bridges. We conclude by extending the definition of $\rho$ to the time interval $[0, T]$ for any $T > 0$ and prove a Central Limit Theorem for the case of two independent Ornstein-Uhlenbeck processes. Indeed, this closes all previously open problems raised in Section 3.3 of \cite{ernst2017yule}. \\

\indent We proceed with some necessary notation. Let $(X_t)_{0 \leq t \leq T}$ be some process with values in $\R^d$, defined over a fixed time interval $[0,T]$. Define random variables
\begin{equation}
\bar{X}  \define T^{-1}\int_0^T X_s \; ds, \qquad Y \define \int_0^T (X_s - \bar{X})
(X_s - \bar{X})^T \;ds
\label{eq1}
\end{equation}
with values in $\R^d$ and $\bbM^d$ respectively, where $\bbM^d$ is the space of $d \times d$ real matrices. Let $Y_{ij}$ be the $(i, j)$-th entry of the matrix $Y$. 
In the case $d=2$, when $X$ is a two-dimensional Wiener process\footnote{Scaling 
properties of Brownian motion show that the law of $\rho$ does not depend on the choice of $T>0$.}, Yule's nonsense correlation can be expressed as
\begin{equation}\label{eq:yule1}
\rho \define  \frac{Y_{12}}{\sqrt{Y_{11}}\sqrt{Y_{22}}}.
\end{equation}
 \tcp{The authors of \cite{ernst2017yule} explicitly calculated the second moment of Yule's nonsense correlation. However, their methodology did not successfully extend to explicitly calculating higher order moments of $\rho$.\\
\indent The key vision of \cite{ernst2017yule} was to calculate Laplace transform of the trivariate
object formed of the three quadratic/bilinear forms of $W$, and to express that transform using Fredholm integral equations. In the present work, we rely instead on the characterization of the moment generating function of the random vector $(Y_{11}, Y_{12}, Y_{22})$. }
This approach inherits from an older and well-developed literature, on the laws of quadratic functionals of Brownian motion. There is a fine survey \cite{donati1997some} which presents the state of the subject as it was in 1997. A range of techniques is available to characterize the laws of quadratic functionals of Brownian motion, including:
\begin{enumerate}
\item eigenfunction expansions --- see, for example,   \cite{chan1991indefinite,chan1994polymer,ernst2017yule, fixman1962radius, levy1951wiener, mac1989extension};
\item identifying the covariance of the Gaussian process as the Green function of a symmetrizable Markov process --- see, for example, \cite{chan1994polymer,dynkin1980markov};
\item stochastic Fubini relations --- \tcp{see, for example, \cite{Donati1993,donati1997some}};
\item It\^o's formula --- see \cite{rogers1992quadratic};
\end{enumerate}
\indent The first of these techniques is historically the first; using it to deliver a simple closed-form solution depends on spotting a simpler form for an infinite expansion. The second works well if we can see a Markov process whose Green function is the covariance of the Gaussian process of interest. The third again requires an insight to transform the problem of interest into a simpler equivalent. The fourth, much less often exploited, deals conclusively with settings where the Gaussian process arises as the solution of a linear stochastic differential equation (SDE); this is the approach we use in the present paper.

Section \ref{sec:asymp} studies the asymptotic properties of $\rho$ as $T \rightarrow \infty$. For this discussion, we will write $X(T), Y(T)$ in place of $X, Y$ defined at \eqref{eq1} and $\rho(T)$ in place of $\rho$ defined at \eqref{eq:yule1} to
emphasize their dependence on the time horizon $T$. In the case of Wiener processes, by the property of self-similarity, it is straightforward to show that $\rho(1)$ and $\rho(T)$ have the same distribution. But for Gaussian processes which are not self-similar, $\rho(T)$ will depend on the value of $T$. Section \ref{sec:asymp} investigates this statistic's asymptotic behavior as $T \rightarrow \infty$. The key result is given by Theorem \ref{thm4}, which proves that, in the case of two independent Ornstein-Uhlenbeck processes, $\sqrt{T}\rho(T)$ converges in distribution as $T \rightarrow \infty$ to a zero-mean Gaussian.

\bigbreak

To summarize then, the main contributions of this paper are, in order of importance:
\begin{enumerate}
\item We characterize the distribution of Yule's nonsense correlation through its moments up to order 16. With these moments in hand, we provide the first density approximation to Yule's nonsense correlation. \textit{This closes the final longstanding open question on the distribution of Yule's nonsense correlation from Yule's 1926 paper (\cite{yule1926})}.
\item We develop the most general version Theorem \ref{thm1} of the `mechanical' It\^o-formula approach of \cite{rogers1992quadratic}, and show how it provides straightforward proofs of all the results of this paper.
\item We find the asymptotic behavior of $\rho(T)$ for the case of independent OU processes. \tcp{We shall see that the variance of $\sqrt{T}\rho(T)$ shrinks to zero as the mean reversion parameter tends  towards infinity. This shows that the empirical correlation coefficient may be viewed as ``sensible'' for testing independence of pairs of OU processes.}
\end{enumerate}

\medskip

\section{Quadratic functionals of Gaussian diffusions}\label{QF}
We shall use the notation $\sdp$ for the space of {\em strictly} positive-definite symmetric $d \times d$ matrices, with the canonical ordering $A \geq B$ meaning that $A-B$ is non-negative definite. The main result is the following.

\begin{theorem}\label{thm1}
Suppose that $\sigma:  [0, T] \mapsto \bbM^d$ is a bounded measurable function, and that  $X$ solves\footnote{For notational simplicity, we will often omit the independent variable $t$. 
}
\begin{equation}\label{eq:linear.sde0}
d X =  \sigma d W, 
\end{equation}
 where $W$ is $d$-dimensional Brownian motion.
 We write $\Sigma = \sigma \sigma^\top$. 
 
 Suppose that $Q: [0,T] \rightarrow \sdp$ and $ z: [0,T] \rightarrow \bbR^d$ are bounded measurable functions such that $Q^{-1}$ is also bounded. Define
\begin{align}
\ell \define \;&     \frac{1}{2} X\cdot Q X + z\cdot X, \label{ldef0}\\ 
F(t,x)  \define \;&  E \left[ \exp \left\{ -  \int_t^T \ell(s)ds  - \ell(T) \right\} \,\bigg| \, X(t) = x \right]\label{Fdef0}. 
\end{align} 
Then $F(t, x)$ is given explicitly as 
\begin{equation}\label{eq:guess0}
F(t, x) = \exp \left( -\frac{1}{2} x\cdot V(t) x - b(t)\cdot x - \gamma(t)   \right),
\end{equation}
where $V,\, b,\, \gamma$ are obtained as the unique  solutions to the system of ordinary differential equations (ODEs),\footnote{We use an ``overdot" to denote the derivative with respect to $t$.}
\begin{align} 
\dot{V} =\;&  V \Sigma V -  Q ,   \label{eq:ode.Q0} \\
\dot{b} =\;&  V \Sigma b  - z    , \label{eq:ode.b0}  \\
2 \dot{\gamma} =\;&  b^\top \Sigma  b  -\trace (V \Sigma), \label{eq:ode.psi0}
\end{align}
subject to the boundary conditions $V(T) = Q(T), \, b(T) = z(T), \, \gamma(T) = 0$. 
\end{theorem}

\begin{proof}
(i) Notice that $\ell$ is bounded below by $-\half z \cdot Q^{-1} z$, which by hypothesis is bounded below by some constant, therefore $F$ defined by \eqref{Fdef0} is bounded.
\\

\noindent
(ii) The ODE \eqref{eq:ode.Q0} has a unique solution up to possible explosion, as the coefficients are locally Lipschitz. We claim that this solution remains positive-definite for $t \leq T$. Since $Q(T) \in \sdp$, it has to be that there exists some $\varepsilon>0$ such that $V(t) \in \sdp$ for all $t \in [T-\varepsilon,T]$. If $V$ does not remain positive definite, then there exists some non-zero $w \in \R^d$ and a greatest $t^* \leq T-\varepsilon < T$ such that $w\cdot V(t^*) w \leq 0$. But we see from \eqref{eq:ode.Q0} that $w \cdot \dot{V}(t^*) w\leq-  w \cdot Q(t^*) w <0$, contradicting the definition of $t^*$. Hence $V$ remains positive-definite all the way back to possible explosion.  However, we have that
\begin{equation*}
V(t) = Q(T) + \int_t^T \{ \,Q(s) -V(s) \Sigma(s) V(s)  \,\} \; ds
\leq Q(T) + \int_t^T Q(s) \; ds.
\end{equation*}
So by hypothesis $V$ is bounded above and no explosion happens. Since $V$ is continuous on $[0,T]$ and positive-definite everywhere, it follows that $V$ is {\em uniformly} positive-definite on $[0,T]$, that is, $V^{-1}$ remains bounded.

It now follows easily that $b$ and $\gamma$ defined by \eqref{eq:ode.b0} and \eqref{eq:ode.psi0} are unique, continuous and bounded. \\

\noindent
(iii) Now define the process
\begin{equation}
Z_t = \half X_t \cdot V_t X_t + b_t \cdot X_t + \gamma_t,
\end{equation}
and develop
\begin{eqnarray*}
dZ_t &=& ( V_t X_t+b_t, \sigma_t \, dW_t \, dt)
+\half \trace(V_t \Sigma_t) dt
+ \bigl\lbrace \half X_t \cdot\dot{V_t}X_t+ \dot{b_t} X_t + \dot{\gamma_t}\bigr\rbrace
dt,
\\
d \langle Z \rangle_t &=& ( V_t X_t+b_t) \cdot \Sigma_t ( V_t X_t+b_t) dt.
\end{eqnarray*}
Now consider the process 
\begin{equation}
M_t = \exp \biggl( \; -\half \int_0^t  \ell(s) \; ds
 - Z_t \biggr).
\label{Mdef}
\end{equation}
Notice that $M$ is bounded, because $\ell$ is bounded below, and so is $Z$ since we have proved that $V^{-1}$, $b$ and $\gamma$ are all bounded on $[0,T]$.
Developing $M$ using It\^o's formula, with the symbol $\doteq$ denoting that the two sides of the equation differ by a local martingale and omitting explicit appearance of the time parameter, we obtain
\begin{eqnarray*}
\frac{dM_t}{M_t}  &=& -dZ + \half d\langle Z \rangle - \half X\cdot Q X dt
-z \cdot X dt
\\
&\doteq & \bigl\lbrace\;
 -\half \trace(V \Sigma)
- \half X \cdot\dot{V}X- \dot{b} X - \dot{\gamma} +
\\
&& \qquad
+\half ( V X+b) \cdot \Sigma ( V X+b)
-\half X \cdot Q X - z \cdot X
\;\bigr\rbrace \; dt
\\
&=& 0
\end{eqnarray*}
because of \eqref{eq:ode.Q0}, \eqref{eq:ode.b0} and \eqref{eq:ode.psi0}.  Thus $M$ is a local martingale, which is also bounded on $[0,T]$ so $M$ is a bounded martingale, and the result follows.

\end{proof}

Theorem \ref{thm1} extends easily to the situation where $X$ is the solution of a linear SDE.

\begin{theorem}\label{thm2}
Suppose that $\sigma, \; B:  [0, T] \mapsto \bbM^d$ and $\delta : [0,T] \mapsto \R^d$ are  bounded measurable functions, and that  $X$ solves
\begin{equation}\label{eq:linear.sde}
d X =  \sigma\,  d W  + (BX +\delta) \, dt, 
\end{equation}
Suppose that $Q: [0,T] \rightarrow \sdp$ and $ z: [0,T] \rightarrow \bbR^d$ are bounded measurable functions such that $Q^{-1}$ is also bounded, and suppose that $\ell$ and $F$ are defined as before at \eqref{ldef0}, \eqref{Fdef0}.

Then $F(t, x)$ is given explicitly as 
\begin{equation}\label{eq:guess}
F(t, x) = \exp \left( -\frac{1}{2} x\cdot V(t) x - b(t)\cdot x - \gamma(t)   \right),
\end{equation}
where $V,\, b,\, \gamma$ are obtained as the unique  solutions to the system of ordinary differential equations (ODEs),\footnote{We use an ``overdot" to denote the derivative with respect to $t$.}
\begin{align} 
\dot{V} =\;&  V \Sigma V - ( V B + B^\top V ) - Q ,   \label{eq:ode.Q} \\
\dot{b} =\;& ( V \Sigma - B^\top) b - V \delta - z    , \label{eq:ode.b}  \\
2 \dot{\gamma} =\;&  b^\top \Sigma  b  -\trace (V \Sigma) - \delta^\top b , \label{eq:ode.psi}
\end{align}
subject to the boundary conditions $V(T) = Q(T), \, b(T) = z(T), \, \gamma(T) = 0$. 
\end{theorem}

\begin{proof}
The coefficients of the SDE \eqref{eq:linear.sde} are globally Lipschitz, so it is a standard result (see, for example, \cite{RW} Theorem V.11.2) that the SDE has a unique strong solution.  If we now set
\begin{equation}
\tilde{X_t} =  A_t X_t + c_t,
\label{tildeX}
\end{equation}
where $A$ and $c$ solve
\begin{eqnarray}
\dot{A}_t + A_t B_t &=& 0, \qquad A(0) = I,
\label{Aode}
\\
\dot{c}_t + A_t \delta_t &=& 0, \qquad c(0) = 0,
\label{code}
\end{eqnarray}
then a few simple calculations show that 
\begin{equation*}
d\tilde{X} = A \sigma \, dW
\end{equation*}
and Theorem \ref{thm1} applies. The equations \eqref{eq:ode.Q}, \eqref{eq:ode.b} and \eqref{eq:ode.psi} are easily checked to be the analogs of \eqref{eq:ode.Q0}, \eqref{eq:ode.b0} and \eqref{eq:ode.psi0} respectively.

\end{proof}

\begin{remark}
We will want to apply Theorem \ref{thm1} to situations where $Q(T) = 0$. This is a simple limiting case of the problem where we take $Q(T) = \varepsilon I$ and let $\varepsilon \downarrow 0$. In a little more detail, we let $V^\varepsilon, \, b^\varepsilon,\, \gamma^\varepsilon$ denote the solution to \eqref{eq:ode.Q}-\eqref{eq:ode.psi} with boundary condition $Q(T) = \varepsilon I$, and we write
\begin{equation}
q^\varepsilon_t : x \mapsto \half x \cdot V^\varepsilon(t)x + b^\varepsilon(t) \cdot x + \gamma^\varepsilon(t),
\end{equation}
for the quadratic form $-\log F(t,x)$. Evidently $q^\varepsilon_t(x)$ is decreasing in $\varepsilon$ for each $x$ and $t$, and from this it follows easily that limits of $V^\varepsilon(t), \, b^\varepsilon(t), \, \gamma^\varepsilon(t)$ exist for each $t$ and determine $F$ for the limit case when $Q(T)=0$.

\end{remark}

\begin{remark}\label{rmk:general}
Theorem~\ref{thm2} is a special case of the Feynman-Kac formula; the fact that the process $M$ defined in \eqref{Mdef} is a martingale is equivalent to the Feynman-Kac formula, and is valid for any additive functional $\ell$ of the diffusion $X$. However, without the special linear form of the SDE for $X$ and the quadratic form of the additive functional $\ell$ it is rare that any explicit solution can be found for $F$.
\end{remark}

\begin{remark}\label{rmk:trans}
If $\sigma$ is constant, we may assume that $\sigma = I$, the identity matrix. 
To see this, let $\hat{X} = \Sigma^{-1/2} X$, and note that the diffusion process $\hat{X}$ solves the linear SDE,
$$ d \hat{X} = ( \Sigma^{-1/2} B \Sigma^{1/2}  \hat{X} + \Sigma^{-1/2} \delta) dt + d W. $$
Letting $\hat{Q} = \Sigma^{1/2} Q \Sigma^{1/2}$ and $\hat{z} = \Sigma^{1/2} z$ we obtain
$$ \ell = \frac{1}{2} X^\top Q X + z^\top X =  \frac{1}{2} {\hat{X}}^\top \hat{Q} \hat{X} + \hat{z}^\top \hat{X},$$
and thus we can work with the process $\hat{X}$ instead of $X$.  However, it seems simpler to provide the full form of the solution for the SDE~\eqref{eq:linear.sde} rather than a reduced form which then requires a translation back to the original problem. 
\end{remark}

\begin{remark}
Although Theorem \ref{thm2} deals with the general case where $Q, z$ are measurable functions, in the remainder of this paper we only need invoke Theorem \ref{thm2} for the special case in which $Q$ and $z$ are constants. For this reason, we will sometimes use the alternative expanded notation
\begin{equation}
F(t,x) \define F(t,x; Q,z)
\label{longFdef}
\end{equation}
when we want to make explicit the dependence of $F$ on the coefficients $Q$ and $z$ appearing in $\ell$.
\end{remark}

\section{Computing the moments of $\rho$}\label{sec:quad}
Henceforth, we deal exclusively with cases where 
\[
d=2.
\]
Recall the definition \eqref{eq1} of the $2 \times 2 $ random matrix $Y$.
Let  $\phi$ be the moment generating function of the joint distribution of $(Y_{11}, Y_{12}, Y_{22})$, which can be expressed using quadratic functionals of $X$ as
\begin{equation}\label{eq:def.exp}
\begin{aligned}
\phi(S ) \define \;&   E  \left[ \exp \left\{  -\half  (s_{11} Y_{11} +  2 s_{12} Y_{12} + s_{22} Y_{22} ) \right\} \right] \\
=\;&   
E  \left[  \exp \left\{  -\frac{1}{2}\int_0^T (X(u) - \bar{X}) \cdot S (X(u)  - \bar{X})  du      \right\} \right]. 
\end{aligned}
\end{equation} 
Here,  $S$ is a $2 \times 2$ positive-definite symmetric matrix with entries denoted by $s_{ij} \; (i,j=1,2).$  As we shall show in the following proposition, the function $\phi$ is all we shall need to evaluate the moments of $\rho$.

\begin{proposition}\label{prop:mgf}
Let $\rho$ be as given in~\eqref{eq:yule1} and $\phi(s_{11}, s_{12}, s_{22}) = \phi(S)$ be as given in~\eqref{eq:def.exp}. For $k = 0, 1, 2, \dots$, we have 
\begin{equation}\label{eq:moments}
E  \rho^{k}  
= \dfrac{(-1)^k}{2^{k} \Gamma(k/2)^2 } \int_0^\infty \int_0^\infty 
s_{11}^{k/2 - 1} s_{22}^{k/2 - 1} \dfrac{\partial^{k} \phi }{ \partial s_{12}^{k} } (s_{11}, 0, s_{22}) \, d s_{11} \, d s_{22}. 
\end{equation}
\end{proposition}

\begin{proof}
It is well known that the moments of a random variable can be obtained by differentiating the moment generating function, given it exists~\cite{billingsley2008probability}. Now note that for any fixed nonnegative $s_{11}, s_{22}$, there exists $\epsilon > 0$ such that $S = [s_{ij}]$ is positive semi-definite for any $s_{12} \in [-\epsilon, \epsilon]$ and thus $\phi(s_{11}, s_{12}, s_{22}) \leq 1$. Hence, the partial derivative with respect to $s_{12}$ exists at $s_{12} = 0$. Applying Fubini's Theorem we obtain
\begin{align*}
(-1)^k \dfrac{\partial^{k} \phi }{ \partial s_{12}^{k} } (s_{11}, 0, s_{22})  
=    E \left[ Y_{12}^{k} \exp\left\{  -\dfrac{1}{2} ( s_{11} Y_{11} + s_{22} Y_{22}  \right\} \right]. 
\end{align*}
Next, recall that by the definition of Gamma function, for any $\alpha > 0$, 
\begin{align*}
y^{- \alpha} = \dfrac{1}{\Gamma(\alpha) } \int_0^\infty  t^{\alpha - 1} e^{- t y} dt
=  \dfrac{1}{2^\alpha \Gamma(\alpha) } \int_0^\infty  s^{\alpha - 1} e^{- s y / 2} ds.
\end{align*}
Since $\rho^k  = Y_{12}^k  \, Y_{11}^{-k/2}  \, Y_{22}^{-k/2} $, we can apply the above formula to obtain~\eqref{eq:moments} (by Tonelli's Theorem, the order of integration can always be exchanged).
\end{proof}

\begin{remark}\label{rmk:mgf}
The idea of using the moment generating function to compute negative moments or moments of the ratio of two random variables has been widely used in the literature. See~\cite{cressie1981moment, jones1987inverse, sawa1972finite}.
\end{remark}

So we see that the distribution of $\rho$ is determined by \eqref{eq:def.exp}, from which moments can in principle be derived using Proposition \ref{prop:mgf}; but we need to get hold of the expression \eqref{eq:def.exp}.  This is where Theorem \ref{thm2} comes in. If $X$ is a solution of a linear SDE \eqref{eq:linear.sde}, starting at $X_0=0$ to fix the discussion, and we set
\[
Q(t) = S, \quad z(t) = a \in \R^2\quad \forall 0 \leq t < T, \qquad Q(T) = 0, \quad z(T) = 0,
\]
then Theorem \ref{thm2} tells us how to compute
\begin{eqnarray}
F(0,0;a) &=& E \biggl[ \exp \biggl\lbrace\; 
  -\int_0^T \{ \half X(u) \cdot S X(u) + a \cdot X(u) \} \; du
\;\biggr\rbrace \biggr]
\label{eq23}
\\
&=& \exp ( -\gamma(0;a) ),
\label{eq24}
\end{eqnarray}
where we have written $F(t,x;a)$ and  $\gamma(0;a)$  to emphasize dependence on $a$. If we now integrate over $a$ with a  $N(0,T^{-1} S)$ distribution the right-hand side of \eqref{eq23} becomes
\begin{equation}
E \biggl[ \exp \biggl\lbrace\; 
  -\int_0^T \half X(u) \cdot S X(u) \; du
  +\half T \bar{X}\cdot S \bar{X}
\;\biggr\rbrace \biggr]  = \phi(S).
\label{eq25}
\end{equation}
The strategy now should be clear. In any particular application, we use Theorem \ref{thm2} to obtain $\gamma(t;a)$ as explicitly as possible, and then we integrate \eqref{eq24} over $a$ to find $\phi(S)$.

\section{Examples}\label{sec4}
In this section we will carry out the program just outlined in four examples, and obtain remarkably explicit expressions for everything we need.

In the first three examples, the two-dimensional diffusion process $X$ has two special properties:
\begin{itemize}
\item[(i)] The law of $(RX_t)_{0 \leq t \leq 1}$ is the same as the law of $(X_t)_{0 \leq t \leq 1}$ for any fixed rotation matrix $R$;
\item[(ii)] The two components of $X$ are independent.
\end{itemize}
Consequently, if we abbreviate $X^1(t) = x_t$, $\bar{x} = \int_0^1 x_s\; ds$, and define
\begin{equation}
\psi(v) = E  \left[  \exp \left\{  -\half \int_0^1 v (x_u - \bar{x})^2   du      \right\} \right],
\label{psidef}
\end{equation}
it follows that the function $\phi(S)$ defined at \eqref{eq:def.exp} simplifies to the product
\begin{equation}
\phi(S)  = \psi(\theta_1^2)\, \psi(\theta_2^2),
\label{phiprod}
\end{equation}
where $\theta_1^2, \, \theta_2^2$ are the eigenvalues of $S$. This observation simplifies the solution of the differential equations \eqref{eq:ode.Q}-\eqref{eq:ode.psi} considerably, reducing everything to a one-dimensional problem.

The final example, that of correlated Brownian motion, reduces to the Brownian example by linear transformation.

\subsection{Brownian motion}\label{ssBM}
For a standard one-dimensional Brownian motion $x(t)$, consider the function $F (t, x; \, \theta^2, z)$ where $\theta \geq 0$ and $z \in \bbR$. By Theorem~\ref{thm2},  the solution has the following form (the subscript ``$\BM$" is Brownian motion)
\begin{equation*} 
F_{\BM}   (t, x; \, \theta^2, z) =  \exp\left\{  -\half V  x^2 - b x  - \gamma \right\}, 
\end{equation*}
which leads to the following system of ordinary differential equations 
\begin{align*}
  \dot{V} - V^2 + \theta^2 &= 0 ,   \\
 \dot{b}  - V b + z &= 0,   \\
  2\dot{\gamma} -  b^2 + V &= 0.   
\end{align*}
Using the boundary condition $V(T) = 0$, we obtain 
$$V(t) =  \theta \tanh   \theta \tau,$$ 
where $\tau = T - t$.  Using the condition $b(T) = 0$, one can show that the solution for $b$ is 
\begin{align*}
b(t) =  \frac{z}{\theta^2} V(t) =   \frac{z}{\theta } \tanh   \theta \tau \, . 
\end{align*}
Solving the third ODE, we obtain  
\begin{align*}
2 \gamma(t) = \log \cosh \theta \tau + \frac{z^2}{\theta^3} \left( - \theta \tau + \tanh \theta \tau    \right), 
\end{align*}
and thus
\begin{align*}
F(0,0;\theta^2, z) = \exp \left\{  - \frac{z^2}{2 \theta^3} \left( - \theta T + \tanh \theta T    \right) -  \frac{1}{2}\log \cosh \theta T  \right\}. 
\end{align*}
As at \eqref{eq25}, we now mix this expression over $z \sim N(0,\theta^2/T)$ to discover that in this example the function $\psi$ (defined at \eqref{psidef}) takes the simple explicit form
\begin{equation}\label{eq:bm.sol1}
\psi_\BM(\theta^2) = \left(  \dfrac{ \theta T }{\sinh\theta T  } \right)^{1/2},
\end{equation}
which (after appropriate scaling), is identical to the relation in the third display on p.577 of \cite{Donati1993}, which is in fact a special case of the result (3.10) on p. 251 of \cite{chan1994polymer}. 

From \eqref{phiprod} therefore, the moment generating function $\phi(S)$ is given by  
\begin{equation}\label{eq:bm.sol}
\phi_\BM(S) = \left(  \dfrac{  \theta_1 \theta_2 T^2 }{\sinh \theta_1 T \sinh \theta_2 T } \right)^{1/2},   
\end{equation}
where $\theta_1^2, \theta_2^2$ are the eigenvalues of $S$. These eigenvalues are given in terms of the entries of $S$ as
\begin{equation}\label{eq:theta}
 \theta_i^2 =    \dfrac{1}{2}\left( s_{11} + s_{22} \pm \sqrt{ (s_{11}-s_{22})^2 + 4s_{12}^2 } \right), 
\end{equation}
where $s_{ij}$ is the $(i, j)$-th entry of $S$. We note that the formula in \eqref{eq:bm.sol} is given in more generality in formula (3.b) on page 578 of \cite{Donati1993}.

Consider $E (\rho^k )$ for $k = 0, 1, 2, \dots$. Note that for any $k$, the expectation always exists since  $\rho \in [-1, 1]$. Further, all the odd moments, i.e. $E (\rho^{2k + 1} )$, are zero by symmetry. 
To compute an even moment of $\rho$, we apply formula~\eqref{eq:moments}. For example,  consider the second moment.  Straightforward but tedious calculations yield  
\begin{align*}
E \rho^2 = \int_0^\infty \int_0^{v}
\dfrac{uv \sqrt{ uv } }{ ( v^2  - u^2 )\sqrt{\sinh u \, \sinh v}} \left( \dfrac{1}{u \tanh u}  - \dfrac{1}{v \tanh v}  
- \dfrac{1}{u^2}  + \dfrac{1}{v^2}  \right) \, d u \, d v, 
\end{align*}
where we have applied a change of variables, $u = \sqrt{s_{11}}, v = \sqrt{s_{22}}$. Note that this is exactly the same as the formula provided in \cite[Proposition 3.4]{ernst2017yule}. 

For higher-order moments, the calculation of $\partial^{k} \phi / \partial s_{12}^k$ is extremely laborious. We use \texttt{Mathematica} to perform symbolic high-order differentiation and then the two-dimensional numerical integration. The numerical results are summarized in Table~\ref{table:bm}. The choice of $T$ is irrelevant since the distribution of $\rho(T)$ does not depend on $T$.

\begin{table}[t!]
\begin{center}
 \begin{tabular}{cccccc}
 \toprule 
$k$ & 2 & 4 & 6 & 8  \\
$E \rho^{k}$ & 0.240522  & 0.109177 & 0.060862  & 0.037788  \\
\midrule 
$k$ & 10 & 12 & 14  & 16 \\
$E \rho^{k}$ & 0.025114 & 0.017504  & 0.012641  & 0.009385 \\
\bottomrule 
 \end{tabular}
 \caption{Numerical values of the moments of Yule's nonsense correlation for two independent Wiener processes ($T = 1$). 
 }\label{table:bm}
\end{center}
\end{table}

We proceed to use the numerical values of $E (\rho^{k})$ to approximate the probability density function of $\rho$, which we denote by $f$. Consider a polynomial approximation
\begin{align*}
\hat{f}_k(\rho) = a_{k,0} + a_{k, 1} \rho + a_{k, 2} \rho^2 + \cdots + a_{k, k} \rho^{k}. 
\end{align*}
The coefficients $(a_{k, 0},  \dots, a_{k, k})$ can be computed by matching the first $k + 1$ moments of $\rho$ (including the zero moment which is always equal to $1$).  
This is also known as the Legendre series expansion of $f$, which minimizes the integrated squared error among all polynomial approximants with degree $k$~\cite{burden1997numerical, provost2005moment}. 
The rate of convergence depends on the modulus of continuity of $f$ (see, for example, \cite{saxena1967expansion,suetin1964representation,wang2021much,wang2012convergence}); the theoretical properties of the latter are difficult to investigate via the moment generating function $\phi$. Below we perform some numerical experiments to show that this polynomial approximation strategy indeed provides an efficient solution to calculating the distribution of $\rho$.

\begin{figure}[t!]
\begin{center}
\includegraphics[width=0.46\linewidth]{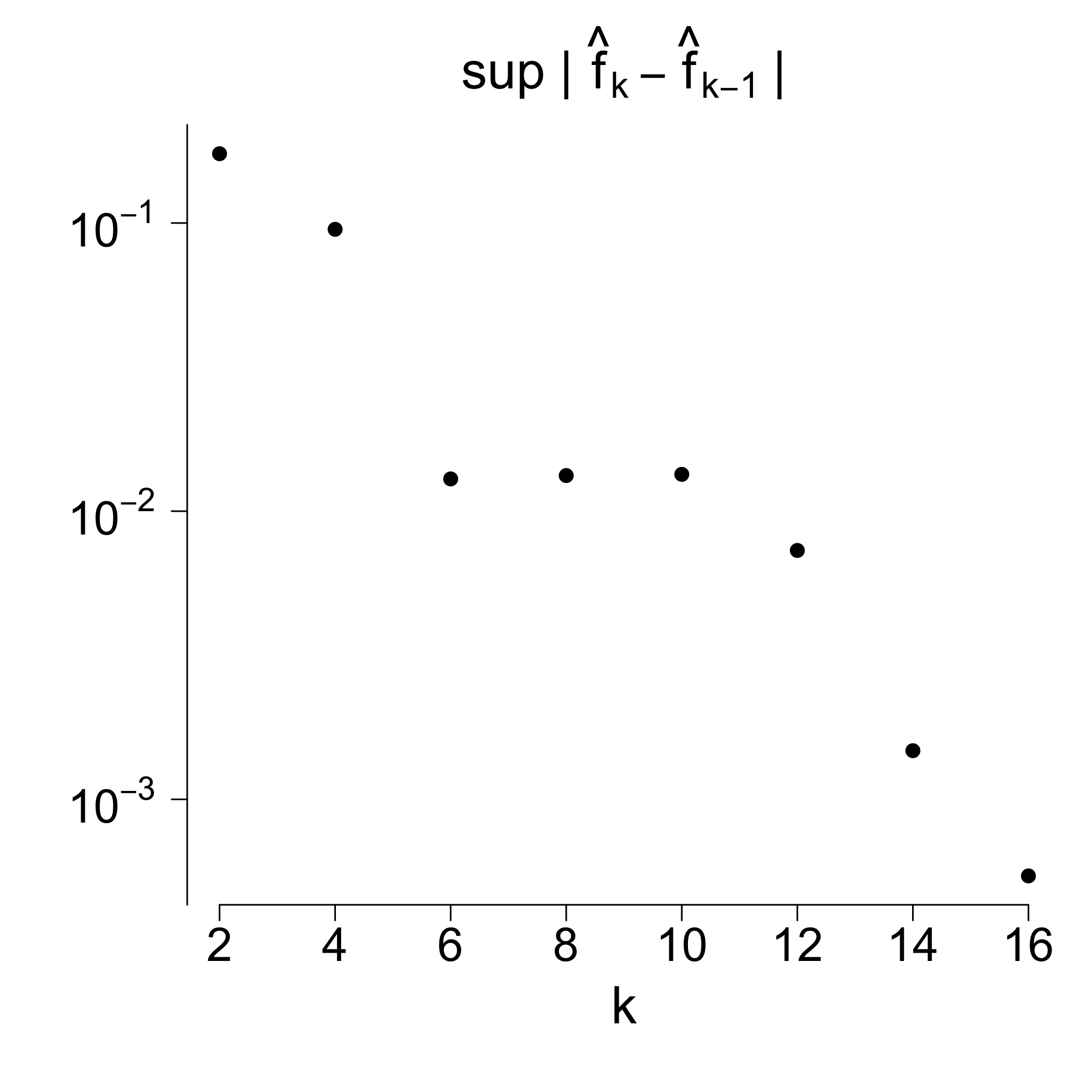}  \hspace{0.6cm}
\includegraphics[width=0.46\linewidth]{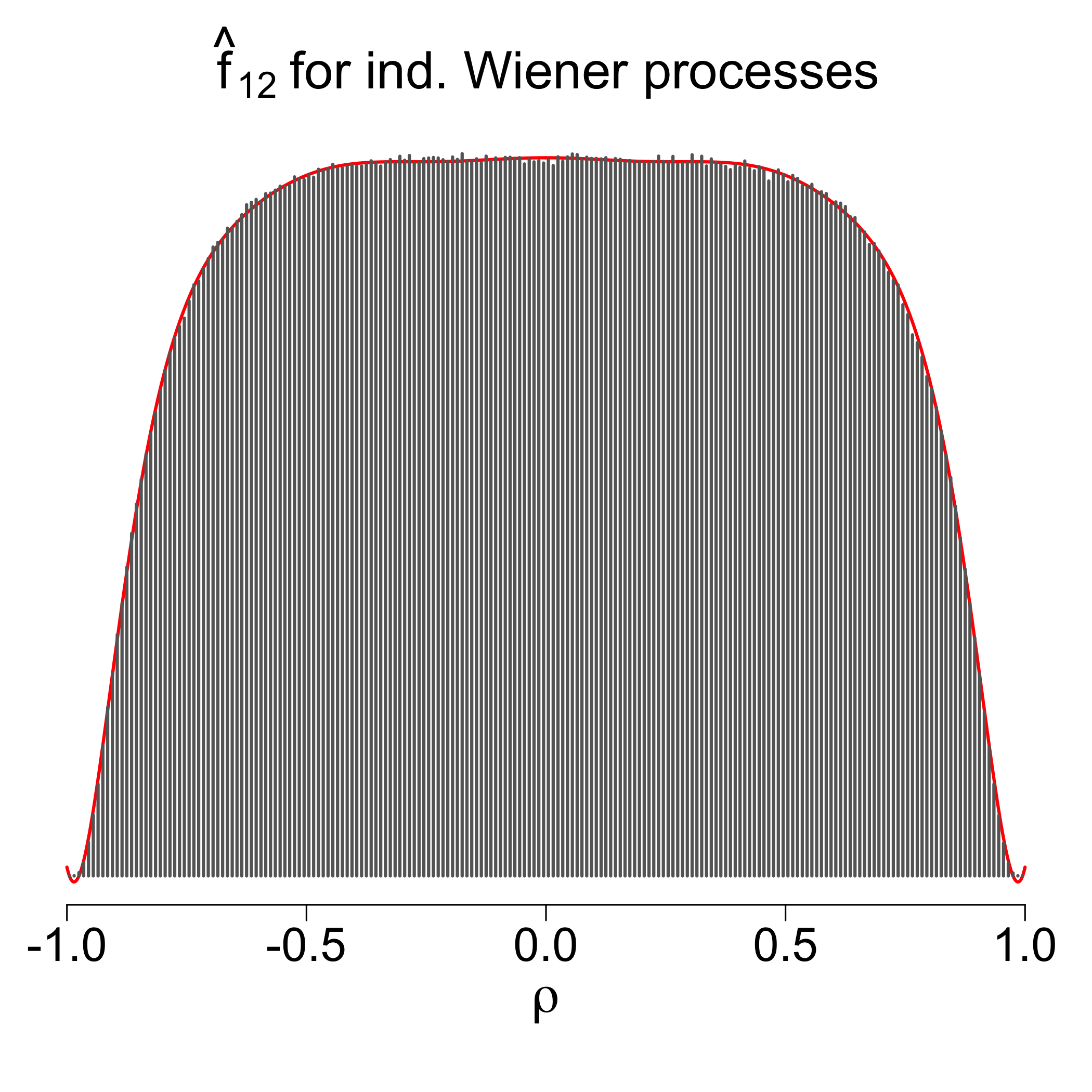} 
\caption{Moment-based polynomial approximations for the density function of $\rho$ for two independent Wiener processes. $\hat{f}_k$ denotes the approximant with degree $k$. In the right-panel, the red curve corresponds to $\hat{f}_{12}$. The gray bars represent the empirical frequencies from $10^7$ Monte Carlo simulations. }\label{fig1}
\end{center}
\end{figure}

Recalling that $E (\rho^k ) = 0$ for odd $k$, we have, for for $k=0, 2, 4,  \dots$, that $\hat{f}_{k}= \hat{f}_{k+1}$. It thus suffices to consider $\hat{f}_{k}$ , for $k=0, 2, 4, \dots$. To determine whether the Legendre series expansion has ``converged'', a commonly used diagnostic is the quantity $\hat{\varepsilon}_k =  \sup_{\rho \in (-1, 1)} |\hat{f}_k(\rho) - \hat{f}_{k-1}(\rho) |$, which needs to be sufficiently small for the algorithm to stop. 
In the left panel of Figure~\ref{fig1}, we show how $\hat{\varepsilon}_k$ changes with $k$ for two independent Wiener processes, from which we see that $\hat{\varepsilon}_k$ quickly tends to zero as $k$ increases. The plot also suggests that $\hat{f}_4$ may give a reasonably good approximation to $f$ since $\hat{\varepsilon}_6 = |\hat{f}_6 - \hat{f}_4 |_\infty$ is just about $0.01$. 
Indeed, we observe that the overall shape of $\hat{f}_4$ is very similar to that of $\hat{f}_{12}$, but $\hat{f}_4$ is more rough. 
Below we give the expressions for $\hat{f}_4, \hat{f}_6$ and $\hat{f}_8$, which constitute \textit{the first density approximations} to Yule's nonsense correlation, and \tcp{thereby \textit{solves the second (and final) of the two longstanding open questions raised by Yule's 1926 paper (\cite{yule1926})!}}
\begin{align*}
\hat{f}_{4}(\rho) =\;& 0.59081 + 0.31001 \rho^2 - 0.97075 \rho^4, \\
\hat{f}_{6}(\rho) =\;& 0.60057 + 0.10518 \rho^2 - 0.35627 \rho^4 - 0.45062 \rho^6, \\
\hat{f}_{8}(\rho) =\;& 0.61200  - 0.30638 \rho^2 + 1.9073 \rho^4 - 4.3742 \rho^6 + 2.1019 \rho^8. 
\end{align*} 
Finally, we compare our moment-based polynomial approximation with Monte Carlo estimates. We fix $T=1$ and discretize time using step size $10^{-4}$.
The gray bars in the right panel of Figure~\ref{fig1} give the histogram of $\rho$ from $10^7$ replicates. 
The red curve is the 12th-order approximation $\hat{f}_{12}$, which agrees very well with the empirical distribution. 
Note that the moment-based polynomial approximation is much more efficient and accurate than the empirical density function obtained from sampling; the latter has two sources of errors, one from Monte Carlo sampling and the other from time discretization.  
From the plot, we see that the distribution of $\rho$ is heavily dispersed and frequently large in absolute value, and that the density remains approximately constant for $\rho \in (-0.5, 0.5)$.  

We have only provided the numerical values of $E (\rho^k)$ up to $k =16$. This has been done for two reasons. Firstly, for practical purposes such density approximation, moments of even higher orders are of much less interest. Secondly, the calculations of the derivative $\partial^k \phi / \partial s_{12}^k$ and the double integral in~\eqref{eq:moments} become extremely slow and require massive memory for $k \geq 16$.

\subsection{Ornstein-Uhlenbeck process}\label{ssOU}
Consider a one-dimensional Ornstein-Uhlenbeck (OU) process which starts from $X(0) = 0$ and evolves according to the following stochastic differential equation:
\begin{equation}\label{eq:sde.ou}
dX(t) = -r X(t) dt + dW(t), \quad \quad r \in (0, \infty). 
\end{equation}
By Theorem~\ref{thm2}, the solution has the form
\begin{equation*} 
F_{\OU}   (t, x; \, \theta^2, z) =  \exp\left\{  -\frac{ 1}{2} V  x^2 - b x  - \gamma \right\}, 
\end{equation*}
which can be obtained by solving the following system of ODEs
\begin{align*}
  \dot{V} - 2r V - V^2 + \theta^2 &= 0 ,   \\
 \dot{b}  - (V  + r )b + z &= 0,   \\
  2\dot{\gamma} -  b^2 + V &= 0.   
\end{align*}
 Using $V(T) = 0$, we solve the first equation to obtain
\begin{align*}
V(t) =  \dfrac{ \theta^2   }{ r  + \eta  \coth \eta \tau } \, , 
\end{align*}
where $\eta = \sqrt{r^2 + \theta^2}.$
The second differential equation is first-order linear, so can be solved explicitly; after some straightforward calculations we obtain
\begin{align*}
b(t) =  \dfrac{ z }{  r + \eta  \coth \eta \tau } \left(  1 + \frac{r}{\eta} \tanh \frac{\eta \tau }{2} \right). 
\end{align*}
Finally, solving the last differential equation yields 
\begin{align*}
2\gamma(t) = \dfrac{z^2}{\theta^2} \left\{
\dfrac{  \left(  1 + \frac{r}{\eta} \tanh \frac{\eta \tau }{2} \right)^2 }{  r + \eta  \coth \eta \tau }
- \dfrac{r^2  }{\eta^3} \tanh \frac{\eta \tau }{2} - \dfrac{ \theta^2 \tau }{\eta^2}    \right\} 
- r \tau  + \log \left( \cosh \eta \tau + \frac{r}{\eta} \sinh \eta \tau \right). 
\end{align*}
Mixing over $z$ with a Gaussian law as before, and using $\tanh (x/2) = \coth x - \csch x$,  we obtain  
\begin{align*}
\psi_\OU(\theta^2; r ) = \sqrt{T} e^{rT/2}    \left\{  \dfrac{\theta^2}{\eta^4} [ 2 r ( \cosh \eta T  - 1) + \eta \sinh \eta T  ] + \dfrac{r^2 T}{\eta^3} [ \eta \cosh \eta T  + r \sinh \eta T ]  \right\}^{-1/2}. 
\end{align*}

If we have two independent Ornstein-Uhlenbeck processes $X_1(t), X_2(t)$ which both start at zero and have common mean reversion parameter $r$, one can check that an orthogonal transformation of $X = (X_1, X_2)$ leaves the joint distribution invariant. Indeed, the new two-dimensional process follows exactly the same SDE. Hence, the moment generating function in this case can be computed by
\begin{align*}
\phi_\OU(S; r) = \psi_\OU(\theta^2_1; r ) \psi_\OU(\theta^2_2; r ),
\end{align*}
where  $\theta_1^2, \theta_2^2$ are the eigenvalues of $S$.

\begin{table}[t!]
\begin{center}
\begin{tabular}{cccccccc}
 \toprule 
$r$ & 0.1 & 0.2 & 0.3 & 0.4 & 0.5 & 1 \\ 
$E \rho^2$ & 0.23209  & 0.22438 & 0.21734  & 0.21091 &  0.20504  & 0.18231\\
\midrule
$r$ & 2 & 5 & 10  & 20 & 50 & 100  \\
$E \rho^2$ & 0.15583 & 0.11454  & 0.07627  & 0.04404  & 0.01907 & 0.00971 \\
\bottomrule 
 \end{tabular}
\caption{Numerical values of the second moment of the empirical correlation coefficient for two independent Ornstein-Uhlenbeck processes with mean reversion parameter $r$  ($T = 1$). }\label{table:ou}
\end{center}
\end{table}

In Table~\ref{table:ou} above we give the numerical values of $E \rho^2$ for independent Ornstein-Uhlenbeck processes with mean reversion parameter $r$  ($T = 1$). Note that as $r \rightarrow \infty$, the processes converge  to constant zero and thus $E \rho^2$ (the variance of $\rho$) goes to zero. Our numerical results show that $E \rho^2 $ decreases slowly.  

\begin{figure}[t!]
\begin{center}
\includegraphics[width=0.32\linewidth]{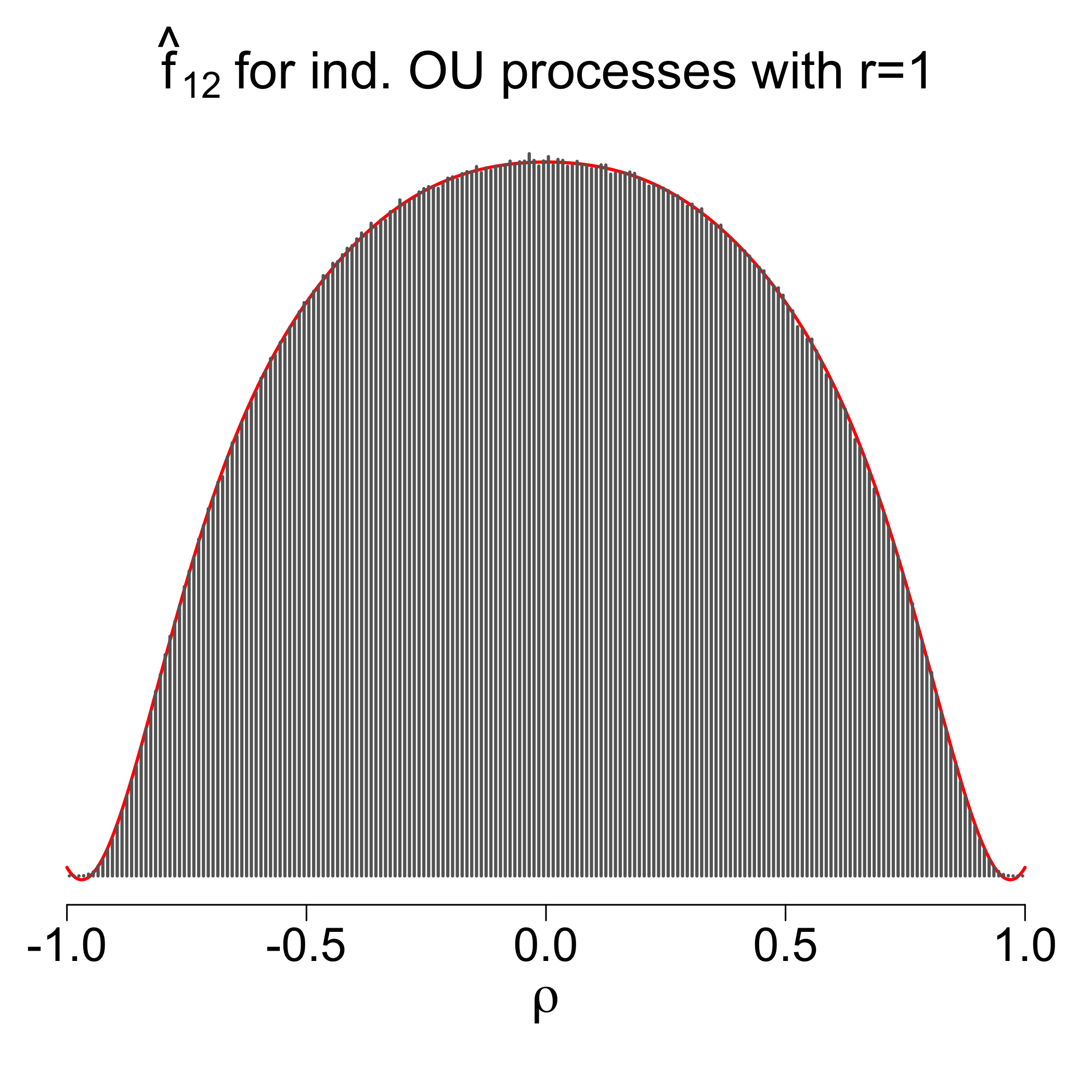}  \hspace{0.1cm}
\includegraphics[width=0.32\linewidth]{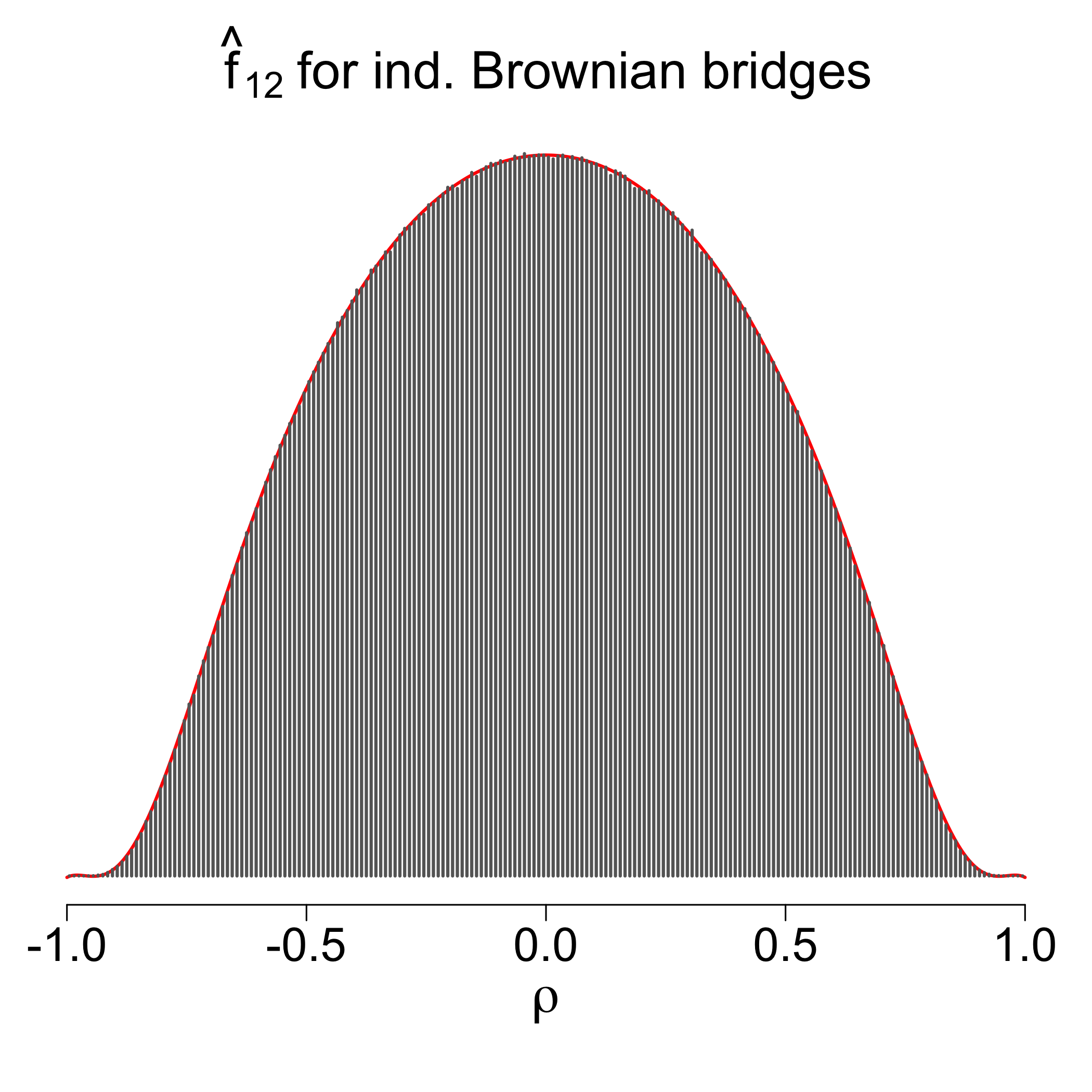}  \hspace{0.1cm}
\includegraphics[width=0.32\linewidth]{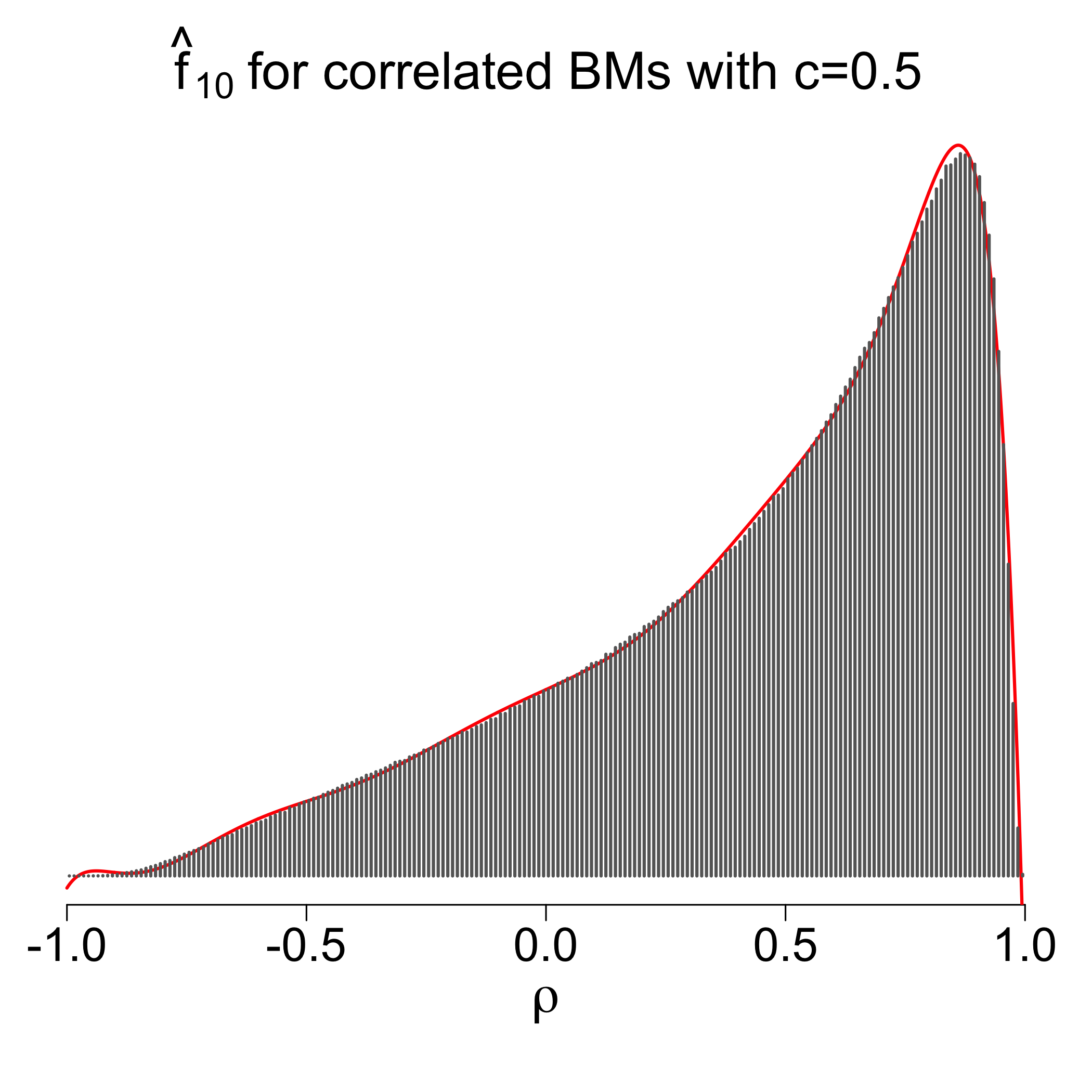} 
\caption{Moment-based polynomial approximations for the density function of $\rho$. 
The gray bars represent the empirical frequencies from $10^7$ Monte Carlo simulations. }\label{fig2}
\end{center}
\end{figure}

As in the case of independent Wiener processes, we can use a Legendre series expansion to approximate the density function of $\rho$. For $r = 1$, the result is shown in the first panel of Figure~\ref{fig2}, from which we see that the $12$-th order approximation is very accurate.

\subsection{Brownian bridge}\label{ssBB}
For a more complicated example, consider a standard  Brownian bridge (denoted by ``$\BB$") which satisfies $X(0) = X(1) = 1$. In this case, we must fix $T = 1$ and let $\tau = 1 - t$. The dynamics of $X(t)$ can be described by
(see for example \cite{RW} Theorem IV.40.3)
\begin{equation*}\label{eq:bb.sde}
d X(t) = -\dfrac{X(t)}{1 - t} dt + d W(t). 
\end{equation*}
Though this SDE has the  linear form, the drift coefficient  $-(1-t)^{-1}$ explodes at $t = 1$. Hence, it does not satisfy the conditions required in Theorem~\ref{thm2}. However, the singularity can easily be isolated, by freezing everything at $t = 1 - \varepsilon$ and applying Theorem \ref{thm2} to that; we can then let $\varepsilon \downarrow 0$ and we find the instances of the ODEs \eqref{eq:ode.Q}-\eqref{eq:ode.psi} to be
\begin{align*}
  \dot{V} - 2V/(1-t) - V^2 + \theta^2 &= 0 ,   \\
 \dot{b}  - [V  + (1-t)^{-1} ] b + z &= 0,   \\
  2\dot{\gamma} -  b^2 + V &= 0.   
\end{align*}
Solving the first differential equation with $\lim_{t \rightarrow 1} V(t) = 0$ yields 
\begin{align*}
V(t) = \dfrac{  \theta \tau \cosh \theta \tau  - \sinh \theta \tau  }{\tau  \sinh \theta \tau}  \, . 
\end{align*}
One can check that $\lim_{t \rightarrow 1} \dot{V}(t) = -\theta^2/3$. 
Similarly, the solution to the second ODE is given by  
\begin{align*}
b(t) = \dfrac{z  (\cosh \theta \tau - 1)}{\theta \sinh \theta \tau } \, ; 
\end{align*}
Though at first sight this might appear to have a singularity at $\tau=0$ it is in fact analytic.
The solution to the third differential equation is given by 
\begin{align*}
2 \gamma(t) = \dfrac{z^2}{\theta^2} \left( \frac{2 (\cosh \theta \tau - 1)}{\theta \sinh \theta \tau }  - \tau \right)+ \log \dfrac{\sinh \theta \tau }{\theta \tau }.
\end{align*}
One can also check that $\lim_{t \rightarrow 1} \gamma(t) = \lim_{t \rightarrow 1} \dot{\gamma}(t) = 0$. 
Using this, we have
\begin{equation*} 
F_{\BB}   (t, x; \, \theta^2, z) =  \exp\left\{  -\frac{ 1}{2} V  x^2 - b x  - \gamma \right\}.
\end{equation*}
Hence
\begin{eqnarray*}
F_{\BB}   (0,0; \, \theta^2, z) &=&  \exp\left\{   - \gamma(0) \right\}
\\
&=& \sqrt{\frac{\theta}{\sinh\theta}} \exp\biggl\lbrace\;
-\frac{z^2}{2 \theta^2}\biggl( \; \frac{2(\cosh \theta -1 )}{\theta \sinh\theta}
-1 \; \biggr)
\;\biggr\rbrace.
\end{eqnarray*}
Mixing over $z \sim N(0,\theta)^2)$ gives the one-dimensional generating function
\begin{align*}
\psi_\BB(\theta^2) = \dfrac{\theta }{2 \sinh (\theta /2)  },
\end{align*}
which matches the formula in the second display on p.577 of \cite{Donati1993}.
As in the case case of Ornstein-Uhlenbeck processes, the moment generating function is
$\phi_\BB(S) = \psi_\BB(\theta^2_1  ) \psi_\BB (\theta^2_2).$

\begin{table}[t!]
\begin{center}
 \begin{tabular}{cccccc}
 \toprule 
$k$ & 2 & 4 & 6 & 8 & 10 \\
$E \rho^{k}$ & 0.149001  & 0.047864  & 0.0201829  & 0.009876 & 0.005321 \\
\bottomrule 
 \end{tabular}
 \caption{Numerical values of the moments of the empirical correlation coefficient for two independent Brownian bridges. 
 }\label{table:bb}
\end{center}
\end{table}

In Table~\ref{table:bb} we provide the moments of $\rho$ for independent Brownian bridges. Comparing with Table~\ref{table:bm}, we can see that $\rho$ has smaller variance for two Brownian bridges. Intuitively, this is because Brownian bridges are forced to fluctuate around zero more frequently than Brownian motions: a Brownian bridge has to return to zero at $t = 1$ but a Brownian motion is likely to make long excursions away from zero. 
The second panel of  Figure~\ref{fig2} shows the $12$-th order polynomial approximation of the density function of $\rho$.

\subsection{Correlated Brownian motion}\label{ssCBM}
Let $X_1(t), \, X_2(t)$ be two Brownian motions with constant correlation $c$, represented by the following SDE
\begin{equation*}
dX_1(t) = dW_1(t),  \quad \quad  dX_2(t) = c dW_1(t) + \sqrt{1 - c^2} dW_2(t). 
\end{equation*}

To compute the moment generating function $\phi(S)$, we take the approach outlined in Remark~\ref{rmk:trans}. 
Define a matrix $M$ as 
\begin{equation*}\label{eq:def.R}
M = M(c) = \begin{bmatrix}
1  &  0 \\
-c (1 - c^2)^{-1/2}  &  (1 - c^2)^{-1/2}
\end{bmatrix}. 
\end{equation*}
Then the process $MX(t) $ is a two-dimensional Brownian motion with independent coordinates. The inverse of $M$ is
\begin{equation*}
M^{-1} = M^{-1}(c) = \begin{bmatrix}
1  &  0 \\
c  &  \sqrt{1 - c^2}
\end{bmatrix}. 
\end{equation*}
We now transform the problem to the uncorrelated case by  
\begin{align*}
\phi_\CBM(S)  =  \phi_\BM( (M^{-1})^\top S M^{-1}  ), 
\end{align*}
where we use ``$\CBM$" to indicate that $X$ is a correlated two-dimensional Brownian motion. The solution may be expressed as 
\begin{equation}\label{eq:cbm.sol}
\phi_\CBM (S; c ) = \left(   \dfrac{\lambda_1 \lambda_2}{\sinh \lambda_1 \sinh \lambda_2 } \right)^{1/2}, 
\end{equation}
where $\lambda_1^2, \lambda_2^2$ are the eigenvalues of the matrix $(M^{-1})^\top   S M^{-1} $. 
Straightforward calculations yield
\begin{align*}
\lambda_i^2 =\;&  \dfrac{1}{2}\left\{ s_{11} + s_{22} +  2 c s_{12}   
\pm \sqrt{   (s_{11}-s_{22})^2 + 4(c s_{11} + s_{12})(c s_{22} + s_{12}) } \right\}. 
\end{align*}

In Table~\ref{table:cbm} we give the first and second moments of $\rho$ for two-dimensional correlated Brownian motion with correlation coefficient $c$. Observe that $E (\rho)$ is always slightly smaller than $c$ if $c \in (0, 1)$. The variance of $\rho$, computed as $\mathrm{Var}(\rho) = E \rho^2 - ( E \rho )^2$, is decreasing (as $c$ increases) but very slowly. Indeed, the standard deviation of $\rho$ is $0.49$ for $c = 0$,  $0.41$ for $c = 0.5$ and $0.25$ for $c  = 0.8$. 
In Table~\ref{table:cbm2} we give high-order moments of $\rho$ for $c = 0.5$. It is somewhat surprising that $E (\rho^k)$ remains close to $0.1$ even when $k = 10$.

Approximating the density function of $\rho$ is more challenging than in the previous three examples due to the slow decay of high-order moments of $\rho$ and the asymmetry of the density function of $\rho$. 
In the last panel of  Figure~\ref{fig2},  we plot the $10$-th order approximation for $c=0.5$, which agrees with the empirical Monte Carlo estimates well, although there appears to be some slight difference near the mode; the difference may be caused by the slow convergence of $\hat{f}_k$ and/or the time discretization scheme used in simulation.

\begin{table}[t!]
\begin{center}
\begin{tabular}{ccccccc}
 \toprule 
$c$ & 0 & 0.1 & 0.2 & 0.3 & 0.4 \\ 
$E \rho$   &  0  & 0.08873 & 0.17792  & 0.26804 & 0.35963  \\
$E \rho^2$ & 0.24052  & 0.24550 & 0.26061  & 0.28636 &  0.32368 \\
$\Var(\rho)$  & 0.2405  & 0.2376 &  0.2290  &  0.2145  &  0.1943 \\
\midrule 
$c$ & 0.5 & 0.6 & 0.7  & 0.8 & 0.9 \\
$E \rho$  & 0.45338  &  0.55004 & 0.65071 &  0.75698 & 0.87151 \\
$E \rho^2$ & 0.37407 & 0.43986  & 0.52477  & 0.63509  & 0.78298 \\
$\Var(\rho)$  &  0.1685  & 0.1373 &  0.1013  &  0.0621  &  0.0235 \\
\bottomrule 
 \end{tabular}
\caption{Numerical values of the moments of the empirical correlation coefficient for two correlated Brownian motions with correlation coefficient $c$  ($T = 1$). }\label{table:cbm}
\end{center}
\end{table}

\begin{table}[t!]
\begin{center}
 \begin{tabular}{cccccc}
 \toprule 
$k$ & 1 & 2 & 3 & 4 & 5  \\
$E \rho^{k}$ & 0.4534  & 0.3741  & 0.2603  & 0.2221 & 0.17137  \\
\midrule 
$k$ & 6 & 7 & 8 & 9 & 10  \\
$E \rho^{k}$ & 0.1515  & 0.1253  & 0.1109  & 0.0948 & 0.0848  \\
\bottomrule 
 \end{tabular}
\caption{Numerical values of higher moments of the empirical correlation coefficient for two correlated Brownian motions with correlation coefficient $c=0.5$  and $T = 1$. }\label{table:cbm2}
\end{center}
\end{table}

\section{Asymptotics of $\rho(T)$ as $T \rightarrow \infty$} \label{sec:asymp}
We now extend the definition of the empirical correlation coefficient to the time interval $[0, T]$ for any $T > 0$. In this context, the empirical correlation coefficient may be written as
\begin{equation}\label{eq:yule}
\rho(T) \define  \dfrac{ Y_{12}(T) }{ \sqrt{ Y_{11}(T) Y_{22}(T) } }, 
\end{equation}
where the random variables $Y_{ij}(T) \, (i, j = 1, 2)$  are defined as  
\begin{equation}\label{eq:def.Y1}
Y_{ij}(T)   \define   \int_0^T  X_i(u) X_j(u) du -  T \, \bar{X}_i    \bar{X}_j   ,  \quad 
\bar{X}_i  \define \dfrac{1}{T} \int_0^T X_i(u) du. 
\end{equation}
The random variable $\bar{X}_i$ is the time average of the process $X_i$. \\
\indent The fundamental reason that the statistic $\rho(T)$ has been called ``nonsense correlation'' is because, in the case of two independent Wiener processes, its asymptotic distribution is heavily dispersed and frequently large in absolute value, leading to high variance (.240522). Further, its asymptotic distribution is very different than that of the nominal $t$-distribution. This begs the following question: might $\rho(T)$ be useful for testing the independence of some \textit{other pair} of Gaussian processes?  In fact, the answer is \textit{yes}; $\rho(T)$ may  be used to test independence of two \textit{Ornstein-Uhlenbeck processes}. We shall prove this claim by first showing a Strong Law result (Theorem \ref{thm4bis}), that for two independent Ornstein-Uhlenbeck processes, $\rho(T)$ converges almost surely to $0$ as $T \rightarrow \infty$. We next prove a Central Limit result (Theorem \ref{thm4}), that $\sqrt{T}\rho(T)$ converges in distribution as $T \rightarrow \infty$ to a zero-mean Gaussian\footnote{Of course, the Strong Law result Theorem \ref{thm4bis}
is not needed to prove the Central Limit result Theorem \ref{thm4}, but as the proof is simple we record it.} with variance that shrinks to zero as the mean reversion parameter tends towards infinity.

\begin{theorem}\label{thm4bis}
For two independent Ornstein-Uhlenbeck processes, $X_1(t)$ and $X_2(t)$, which both follow the SDE~\eqref{eq:sde.ou} with $r > 0$, $\rho(T)$ converges almost surely to zero as $T \rightarrow \infty$. 
\end{theorem}

\begin{proof}
If $X_1(0)$ and $X_2(0)$ are both distributed according to the invariant $N(0,1/2r)$ distribution of the OU process \eqref{eq:sde.ou}, then the bivariate process  $(X_1, X_2 )$ is ergodic, so, by Birkhoff's Ergodic Theorem, time-averages converge almost surely to expectations. Thus (recalling \eqref{eq:yule} and \eqref{eq:def.Y1}) we have
\begin{eqnarray*}
\bar{X_i} &\rightarrow& E[X_i(0)] \qquad= 0 \quad \hbox{\rm a.s. as $T\rightarrow
\infty$, $i=1,2$}
\\
T^{-1}Y_{12}(T)&\rightarrow& E[X_1(0)X_2(0)] = 0  \quad \hbox{\rm a.s. as $T\rightarrow
\infty$}
\\
T^{-1}Y_{ii}(T)& \rightarrow& E[X_i(0)^2] \qquad = (2r)^{-1}\quad
\hbox{\rm a.s. as $T\rightarrow
\infty$, $i=1,2$}
\end{eqnarray*}
Dividing the numerator and denominator of $\rho(T)$ defined at \eqref{eq:yule} by $T$, it is immediate that $\rho(T)$ converges almost surely to 0 if the initial distribution is the invariant distribution.

 If the initial distribution is something else, then we still have these results by coupling with an independent stationary copy of the OU process - see \cite{RW} Theorem V.54.5, which proves that the two diffusions couple in finite time almost surely, so that the long-time averages have the same limits.

\end{proof}

We now prove a central limit theorem for $\rho(T)$ as $T \rightarrow \infty$.

\begin{theorem}\label{thm4}
For two independent Ornstein-Uhlenbeck processes, $X_1(t)$ and $X_2(t)$, which both follow the SDE~\eqref{eq:sde.ou} with $r > 0$, we have that

\beqn
\sqrt{T}\rho(T) \convd N\parens{0, \frac{1}{2r}}.
\eeqn
\end{theorem}

\begin{proof}
Firstly, as we proved in the previous result, we have

\beqn
\begin{aligned}
\frac{Y_{11}(T)}{T}=\;&\frac{1}{T}\int_0^T X_1(s)^2ds-\bar{X}_1(T)^2 \\
\convas\;& E[X_1(0)^2]=\frac{1}{2r}.
\end{aligned}
\eeqn
We now need to obtain weak convergence of

\bneqn \label{eqforweakconv}
\frac{Y_{12}(T)}{\sqrt{T}}=\int_0^T X_1(s)X_2(s)\frac{ds}{\sqrt{T}}-\frac{T\bar{X_1}(T)\bar{X_2}(T)}{\sqrt{T}}.
\eneqn
Let us first consider the second term of the right-hand side of the above equation. For simplicity, assume $X_1(0)=0$ and then

\beqn
X_1(t)=e^{-r t}\int_0^t e^{r s}dW_s\,,
\eeqn
so that
\beqn
\begin{aligned}
T\bar{X_1}(T)=\int_0^T X_1(t)dt=\;&\int_0^Te^{-r t}\int_0^t e^{r s}dW_s \,dt \\
=\;& r^{-1} \int_0^T e^{r s} \parens{-e^{-r T}+e^{-r s}}dW_s\\
=\;& r^{-1} \int_0^T \parens{1-e^{-r (T-s)}}dW_s.
\end{aligned}
\eeqn
Hence
\beqn
E[\bar{X_1}(T)^2]=\frac{1}{r^2T^2}\int_0^T \parens{1-e^{-r (T-s)}}^2 ds \leq \frac{1}{r^2 T},
\eeqn
and so 
\beqn
E\bracks{\parens{\sqrt{T}\bar{X_1}(T)\bar{X_2}(T)}^2}\leq  \frac{1}{r^4 T} \rightarrow 0.
\eeqn
Thus $\sqrt{T}\bar{X_1}(T)\bar{X_2}(T)$ converges in $L^2$ to 0, and so we need now only consider the first term of the right-hand side of equation \eqref{eqforweakconv}. For $\theta \in \mathbb{R}$, let us evaluate the characteristic function by firstly conditioning on $X_2$:

\beqn
\begin{aligned}
&E \exp\left\{ \; \frac{i \theta}{\sqrt{T}}\int_0^T X_1(s)X_2(s)\,ds\right\}\\
=\,&E \exp \left\{ \frac{i \theta}{\sqrt{T}} \int_0^T e^{-r s} X_2(s) \int_0^s e^{r u}dW_1(u)\,ds\right\}\\
=\,&E \exp \left\{ \frac{i \theta}{\sqrt{T}} \int_0^T e^{r u}  \int_u^T  e^{-r s} X_2(s)\,ds\,dW_1(u)  \right\}\\
=\,&E \exp \left\{ -\frac{\theta^2}{2T} \int_0^T \parens{\int_u^T e^{-r(s-u)}X_2(s)ds}^2 du \right\}.
\end{aligned}
\eeqn
Again by the ergodic theorem, we have
\beqn
\frac{1}{T} \int_0^T \parens{\int_u^T e^{-r(s-u)}X_2(s)ds}^2\,du \convas E\bracks{\parens{\int_0^\infty e^{-r s} X_2(s) ds}^2}=\frac{1}{8r^3}.
\eeqn
We thus obtain
\beqn
\frac{1}{\sqrt{T}} \int_0^T X_1(s)X_2(s)ds \convd N\parens{0, \frac{1}{8r^3}},
\eeqn
from which the stated result follows.
\end{proof}


\bigskip

\noindent \textbf{Acknowledgments}  We thank Professor I. Corwin, Professor V. de la Pena,
and Professor F. Viens for many helpful conversations about this work. The first named author acknowledges, with gratitude, the support of Office of Naval Research (ONR) grants N00014-18-1-2192 and N00014-21-1-2672.

\newpage
\bibliographystyle{plain}
\bibliography{yule}

\begin{thebibliography}{10}

\bibitem{Aldrich}
J.~Aldrich.
\newblock Correlations genuine and spurious in {P}earson and {Y}ule.
\newblock {\em Statistical Science}, 10(4):\,364--376, 1995.

\bibitem{billingsley2008probability}
P.~Billingsley.
\newblock {\em Probability and Measure}.
\newblock John Wiley \& Sons, 2008.

\bibitem{burden1997numerical}
R.~L. Burden and J.~D. Faires.
\newblock Numerical analysis.
\newblock {\em Cole, Belmont}, 1997.

\bibitem{chan1991indefinite}
T.~Chan.
\newblock Indefinite quadratic functionals of {G}aussian processes and
  least-action paths.
\newblock {\em Annales de l'IHP Probabilit{\'e}s et Statistiques},
  27(2):239--271, 1991.

\bibitem{chan1994polymer}
T.~Chan, D.~S. Dean, K.~M. Jansons, and L.~C.~G. Rogers.
\newblock On polymer conformations in elongational flows.
\newblock {\em Communications in Mathematical Physics}, 160(2):239--257, 1994.

\bibitem{cressie1981moment}
N.~Cressie, A.~S. Davis, and J.~L. Folks.
\newblock The moment-generating function and negative integer moments.
\newblock {\em The American Statistician}, 35(3):148--150, 1981.

\bibitem{Donati1993}
C.~Donati-Martin and M.~Yor.
\newblock On some examples of quadratic functionals of {B}rownian motion.
\newblock {\em Advances in Applied Probability}, 25:570--584, 1993.

\bibitem{donati1997some}
C.~Donati-Martin and M.~Yor.
\newblock Some {B}rownian functionals and their laws.
\newblock {\em The Annals of Probability}, 25(3):1011--1058, 1997.

\bibitem{dynkin1980markov}
E.~B. Dynkin.
\newblock Markov processes and random fields.
\newblock {\em Bulletin of the American Mathematical Society}, 3(3):975--999,
  1980.

\bibitem{ernst2017yule}
P.~A. Ernst, L.~A. Shepp, and A.~J. Wyner.
\newblock Yule's {``}nonsense correlation'' solved!
\newblock {\em The Annals of Statistics}, 45(4):1789--1809, 2017.

\bibitem{fixman1962radius}
M.~Fixman.
\newblock Radius of gyration of polymer chains.
\newblock {\em The Journal of Chemical Physics}, 36(2):306--310, 1962.

\bibitem{gran}
C.~Granger and D.~Newbold.
\newblock Spurious regression in econometrics.
\newblock {\em Journal of Econometrics}, 2:\,111--120, 1974.

\bibitem{Hendry}
D.~F. Hendry.
\newblock Economic modelling with cointegrated variables: an overview.
\newblock {\em Oxford Bulletin of Economics and Statistics}, 48(3):\,201--212,
  1986.

\bibitem{jones1987inverse}
M.~C. Jones.
\newblock Inverse factorial moments.
\newblock {\em Statistics \& Probability Letters}, 6(1):37--42, 1987.

\bibitem{levy1951wiener}
P.~L{\'e}vy.
\newblock Wiener's random function, and other {L}aplacian random functions.
\newblock In {\em Proceedings of the Second Berkeley Symposium on Mathematical
  Statistics and Probability}. The Regents of the University of California,
  1951.

\bibitem{mac1989extension}
P.~Mac~aonghusa and J.~V. Pule.
\newblock An extension of {L}{\'e}vy's stochastic area formula.
\newblock {\em Stochastics: An International Journal of Probability and
  Stochastic Processes}, 26(4):247--255, 1989.

\bibitem{phil}
P.~C.~B. Phillips.
\newblock Understanding spurious regressions in econometrics.
\newblock {\em Journal of Econometrics}, 33(3):\,311--340, 1986.

\bibitem{phil2}
P.~C.~B. Phillips.
\newblock New tools for understanding spurious regressions.
\newblock {\em Econometrica}, 66(6):1299--1325, 1998.

\bibitem{provost2005moment}
S.~B. Provost.
\newblock Moment-based density approximants.
\newblock {\em Mathematica Journal}, 9(4):727--756, 2005.

\bibitem{rogers1992quadratic}
L.~C.~G. Rogers and Z.~Shi.
\newblock Quadratic functionals of {B}rownian motion, optimal control, and the
  “{C}olditz” example.
\newblock {\em Stochastics: An International Journal of Probability and
  Stochastic Processes}, 41(4):201--218, 1992.

\bibitem{RW}
L.~C.~G. Rogers and David Williams.
\newblock {\em Diffusions, {M}arkov {P}rocesses and {M}artingales: {V}olume 2,
  {I}t{\^o} {C}alculus}.
\newblock Cambridge {U}niversity {P}ress, 2000.

\bibitem{sawa1972finite}
T.~Sawa.
\newblock Finite-sample properties of the $k$-class estimators.
\newblock {\em Econometrica}, 40(4):653, 1972.

\bibitem{saxena1967expansion}
R.~B. Saxena.
\newblock Expansion of continuous differentiable functions in {F}ourier
  {L}egendre series.
\newblock {\em Canadian Journal of Mathematics}, 19:823--827, 1967.

\bibitem{suetin1964representation}
P.~K. Suetin.
\newblock On the representation of continuous and differentiable functions by
  {F}ourier series in {L}egendre polynomials.
\newblock In {\em Doklady Akademii Nauk}, volume 158, pages 1275--1277. Russian
  Academy of Sciences, 1964.

\bibitem{wang2021much}
H.~Wang.
\newblock How much faster does the best polynomial approximation converge than
  {L}egendre projection?
\newblock {\em Numerische Mathematik}, 147(2):481--503, 2021.

\bibitem{wang2012convergence}
H.~Wang and S.~Xiang.
\newblock On the convergence rates of {L}egendre approximation.
\newblock {\em Mathematics of Computation}, 81(278):861--877, 2012.

\bibitem{yule1926}
G.~U. Yule.
\newblock Why do we sometimes get nonsense-correlations between time-series?--a
  study in sampling and the nature of time-series.
\newblock {\em Journal of the Royal Statistical Society}, 89(1):1--63, 1926.

\end{thebibliography}

 \end{document}